\documentclass[11 pt]{article}  

\usepackage[latin1]{inputenc}
\usepackage{amsmath,indentfirst,qsymbols,graphicx,psfrag,amsfonts,stmaryrd,hyperref,color,enumerate}

\def \comp{{\sf Comp}}
\def \1{\textbf{1}}

\def \FD{B}
\def \BSEG{{\sf BSEG}}
\def \SEG{{\sf SEG}}
\def \AC{{\sf AC}}

\def \bar{\overline}
\def \ben{\begin{eqnarray}}
\def \een{\end{eqnarray}}
\def \ba{\begin{align}}
\def \ea{\end{align}}

\def \be{\begin{eqnarray*}}
\def \ee{\end{eqnarray*}}

\def \beq{\begin{equation}}
\def \eq{\end{equation}}

\def \build#1#2#3{\mathrel{\mathop{\kern 0pt#1}\limits_{#2}^{#3}}}

\def \ba{{\bf a}}

\def \sx{{\sf x}}

\def \sz{{\sf z}}

\def \captionn#1{\begin{center}\begin{minipage}{15cm}\sf\caption{\small #1}\end{minipage}\end{center}}

\def \Proof{{\bf Proof. \rm}}

\def \eref#1{(\ref{#1})}

\def \l{\left}

\def \r{\right}
\def \sous#1#2{\mathrel{\mathop{\kern 0pt#1}\limits_{#2}}}
\def \sur#1#2{\mathrel{\mathop{\kern 0pt#1}\limits^{#2}}}

\def \BCH{\partial{\sf Convex Hull}}

\newqsymbol{`B}{\mathcal{B}}
\newqsymbol{`E}{\mathbb{E}}
\newqsymbol{`I}{\mathbb{I}}
\newqsymbol{`N}{\mathbb{N}}
\newqsymbol{`O}{\Omega}
\newqsymbol{`P}{\mathbb{P}}
\newqsymbol{`Q}{\mathbb{Q}}
\newqsymbol{`R}{\mathbb{R}}
\newqsymbol{`W}{\mathbb{W}}
\newqsymbol{`Z}{\mathbb{Z}}
\newqsymbol{`a}{\alpha}
\newqsymbol{`e}{\varepsilon}
\newqsymbol{`o}{\omega}
\newqsymbol{`t}{\tau}
\newqsymbol{`w}{{\cal W}}
\def\cro#1{\llbracket#1\rrbracket}

\def\CP{{\sf CP}}

\begin{document}

\newtheorem{fig}{\hspace{2cm} Figure}
\newtheorem{lem}{Lemma}
\newtheorem{defi}[lem]{Definition}
\newtheorem{pro}[lem]{Proposition}
\newtheorem{theo}[lem]{Theorem}
\newtheorem{cor}[lem]{Corollary}
\newtheorem{note}[lem]{Note}
\newtheorem{conj}{Conjecture}
\newtheorem{Ques}{Question}

\renewcommand{\baselinestretch}{1.2}

 \begin{center}
 \LARGE\sf
Probability that $n$  random points in a disk are in convex position.\\
 {\large \bf Jean-Fran\c{c}ois Marckert}
 \rm \\
 \large{CNRS, LaBRI, Universit\'e Bordeaux \\
  351 cours de la Libération\\
 33405 Talence cedex, France}
 \normalsize
 \end{center}

\begin{abstract} In this paper we give a formula for the probability that $n$ random points chosen under the uniform distribution in a disk are in convex position. While close, the formula is recursive and is totally explicit only for the first values of $n$. 
\end{abstract}~\\ 
{\small \sf Part of this work is supported by ANR blanc PRESAGE (ANR-11-BS02-003)}
 
\section{Introduction}

All the random variables are assumed to be defined on a common probability space $(\Omega,{\cal A},`P)$. The expectation is denoted by $`E$. The plane will be sometimes viewed as $`R^2$ or as $\mathbb{C}$ and we will pass from the real notation (e.g. $(x,y)$) to the complex one $(\rho e^{i \theta})$ without any warning.  For a set $A$ in $`R^2$, $|A|$ denotes the Lebesgue measure of $A$.  
 We denote by $\partial B$ the border of a set $B$. 
For any $n\geq 1$, any $z$, notation $z[n]$ stands for the $n$ tuple $(z_1,\dots,z_n)$ and $z\{n\}$ for the set $\{z_1,\dots,z_n\}$.  
For $H$ be a compact convex domain in $`R^2$ with non empty interior, for any $n\geq 0$, $`P_H^n$ denotes the law of $n$ i.i.d. points $z[n]$ taken under the uniform distribution over $H$. 
A $n$-tuple of points $\sx[n]$ of the plane is said to be in convex position if the $x_i$'s all belong to $\BCH(\sx\{n\})$. Further we define
\[\CP_{n,m}=\big\{\sx[n]: \#\{ i~:~x_i\in \BCH(\sx\{n\})\}=m\big\}\]
the set of $n$ tuples $\sx[n]$ for which exactly $m$ are on the border of $\BCH(\sx\{n\})$.
Hence $\CP_n:=\CP_{n,n}$ is the set of  $n$-tuples of points in convex position. 
Finally, we let 
\ben
P^{n}_H&=&`P_H^n(z[n] \in \CP_n),\\
P^{n,m}_H&=&`P_H^n(z[n] \in \CP_{n,m}).
\een

The aim of the paper is to establish a formula for $P_D^n$, the probability that $n$ i.i.d. random points taken under the uniform distribution in a disk $D$ are in convex position; we will also compute $P_D^{n,m} $ the probability that exactly $m$ points among these $n$ points are on $\BCH(\sz\{n\})$.  
To compute  $P_D^n$ we need and obtain a result more general than the disk case only, result about for what we will call \it bi-pointed segments \rm ($\BSEG$). This will play somehow the role of the bi-pointed triangle (see \eref{eq:bi1}) as studied by B\'ar\'any \& al \cite{BRSZ}, central  also in the approach of Buchta \cite{BC2} (see \eref{eq:bi2}) of the computation of $P^{n,m}_T$ and $P^{n,m}_S$ (where $T$ stands for triangle, and $S$ for square).\par
For $\theta\in[0,2\pi]$, $R> 0$, the \it arc of circle \rm $\AC(\theta,R)$ is defined by 
\[\AC(\theta,R)=\{ R e^{i\nu}, \nu \in[-\theta/2,\theta/2] \}.\]   
\begin{figure}[h]
\centerline{\includegraphics{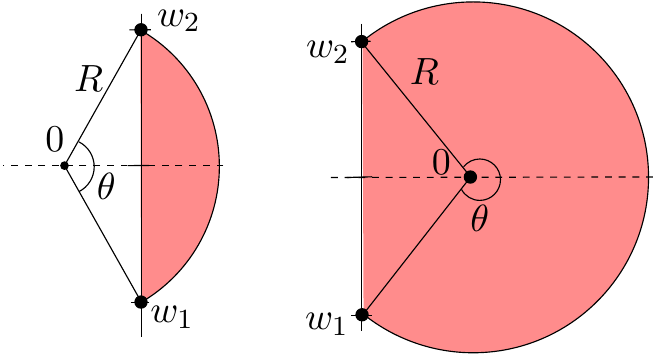}}
\captionn{\label{fig:displ1}Representation of typical $\SEG(\theta,R)$ for $0<\theta<\pi$ and $\pi<\theta<2\pi$}
\end{figure} 
We denote by $\SEG(\theta,R)$ the \it segment \rm (SEG) corresponding to the convex hull of $\AC(\theta,R)$ (see Fig. \ref{fig:displ1}). 
Now consider $w_1(\theta,R)=Re^{-i\theta/2}$ and $w_2(\theta,R)=R e^{i\theta/2}$ the two \it extremities \rm of the \it special border \rm $[w_1(\theta,R),w_2(\theta,R)]$ of $\SEG(\theta,R)$. Let $z_1,\dots, z_n$ be i.i.d. and uniform in $\SEG(\theta,R)$. Set 
\[Z[n,\theta,R]=\{w_1(\theta,R),w_2(\theta,R),z_1,\dots,z_n\},\]
and define the crucial \it bi-pointed segment case \rm ($\BSEG$) function
\ben \label{eq:ere}
\FD_{n,m}(\theta):=`P(Z[n,\theta,R]\in\CP_{n+2,m+2}), ~~~\theta\in(0,2\pi), 1\leq m \leq n.
\een
The value of $R$ has no importance since there exists a dilatation sending $\SEG(\theta,R)$ on $\SEG(\theta,R')$, and dilatations conserve convex bodies and uniform distribution. But it will be useful to have the two parameters $(\theta,R)$ for subsequent computations. 
Again, we write $B_n$ instead of $B_{n,n}$ and below $L_n$ instead of $L_{n,n}$.
Clearly, for any $\theta \in(0,2\pi)$, $\FD_0(\theta)=\FD_1(\theta)=1$. Now for any $n\geq 0$, $\theta\in(0,2\pi)$ define
\ben 
\label{eq:F-L}
L_{n,m}(\theta)&=&\frac{\FD_{n,m}(\theta)(\theta-\sin(\theta))^{n} \sin(\theta/2) }{n!}
\een
Hence 
\ben\label{eq:L01}
L_0(\theta)=\sin(\theta/2),~~~ L_1(\theta)=\sin(\theta/2)(\theta-\sin(\theta)).\een
Notice that $0$ (corresponding to the flat case) as well as  $2\pi$ (the circle case) are excluded from definitions \eref{eq:ere} and \eref{eq:F-L}. 
The main contribution of this paper is the following theorem which allows to compute $P_D^{n,m}$. 
\begin{theo}\label{theo:circle}
\begin{itemize}\itemsep0em 
\item [$(i)$] For any $n\geq 1,$
\[ P_D^n =\lim_{t\to 2\pi\atop{t<2\pi}} \FD_{n-1}(t).\]
\item [$(i')$] For any $n\geq 2$, 
\be
P_D^{n} &=&\frac{(n-2)!}{2^{n-2}\pi^{n-1}} \int_{0}^{2\pi}\sum_{k=0}^{n-2} 
 L_k(\phi)L_{n-2-k}(2\pi-\phi)   d\phi 
\ee
\item [$(ii)$] For any $\theta \in (0,2\pi)$, any $n\geq 1$,
\ben\label{eq:sing}
\frac{L_n(\theta)}{2}= \int_0^{\theta} 
\frac{\sin(\theta/2)^{2n+1}}{\sin(\phi/2)^{2n+1}}  \int_{0}^{\phi} 
 \sum_{k=0}^{n-1}L_{k}(\eta )L_{n-1-k}( \phi-\eta )   d\eta      d\phi.
\een
Analogous results can be obtained for $P_D^{n,m}$:
\item[$(iii)$]  For any $\theta\in(0,2\pi)$, for any $k,l\geq k+1$, $L_{k,l}(\theta)=0$. 
For any $\theta \in (0,2\pi)$, any $n\geq 1$ and $1\leq m \leq n$
\be
\frac{L_{n,m}(\theta)}2&=&\int_{0}^{\theta}\int_{0}^{\phi}\frac{\sin(\theta/2)^{2n+1}}{\sin(\phi/2)^{2n+1}} 
\sum_{n_1+n_2+n_3=n-1\atop{m_1+m_2=m-1}}
 \frac{(\sin(\eta)+\sin(\phi-\eta)-\sin(\phi))^{n_3}}{n_3!}\\
&\times &
L_{n_1,m_1}( \eta)
L_{n_2,m_2}(\phi-\eta)d\eta d\phi
\ee
An alternative form can be given using
\be
\sin(\eta)+\sin(\phi-\eta)-\sin(\phi)&=& 4\sin(\frac{\phi-\eta}2)\sin(\phi/2)\sin(\eta/2).
\ee
\item[$(iii')$] For any $n\geq 2$, $1\leq m \leq n$
\be
P_D^{n,m}&=&\frac{(n-2)!}{2^{n-2}\pi^{n-1}} \int_{0}^{2\pi}\sum_{n_1+n_2=n-1\atop{m_1+m_2=m-1}  }                          L_{n_1,m_1}(\phi)L_{n_2,m_2}(2\pi-\phi)    d\phi.
\ee
\item [$(iv)$] For any $n\geq 1,$
\[ P_D^{n,m} =\lim_{t\to 2\pi\atop{t<2\pi}} \FD_{n-1,m-1}(t).\]

\end{itemize}
\end{theo}
From $(ii)$, one can compute successively the $L_j(\theta)$, and by \eref{eq:F-L}, this allows one to compute the $\FD_n(\theta)$. By $(i)$ it suffices then to take the limit when $\theta\to 2\pi^{-}$ (limit from below). \medskip

Despite important efforts we were not able to find a simpler formula for $\FD_n$ than that presented in the Theorem. 
Nevertheless, explicit computation can be done but close formula for the first $L_j$ given below shows a rapid growth of the complexity of the formula's ($L_{10}$ would need one page to be written down). The effective computation of the first $L_n$ is complex and very few can be computed by hand. In particular the singularity apparent in \eref{eq:sing}  is difficult to handle since the terms in the sum needs to be combined to compensate the singularity. 

In Section \ref{sec:CC} we present an algorithm allowing one to compute the first terms of the sequence (using a computer, or years of time of an efficient human brain).  With this algorithm I computed the 11 first values of $L_n$ which allows the computation $(P_D^n,1\leq n \leq 12)$. $L_0$ and $L_1$ have been given in \eref{eq:L01}; the next ones are
\be
 L_2(\theta)&=&\frac16\sin(\frac{\theta}{2})(3\theta^2+\sin(\theta)^2-16\sin(\frac{\theta}{2})^2)\\
L_3(\theta)&=&
{\frac {1}{54}}\,\sin (\frac{\theta}2)  \left( 2\,  \sin (\frac{\theta}2)  ^{4}\sin ( \theta  ) +9{\theta}^{3}+27\,   \sin (\frac{\theta}2) 
  ^{2}\theta+7    \sin (\frac{\theta}2)^{2}\sin  ( \theta  ) +105 (\sin \left( \theta \right)-\theta)  \right) \\
L_4(\theta)&=&{\frac {\sin (\frac{\theta}2) }{12960}} \Big( 160\,
   \sin (\frac{\theta}2)    ^{6}+48\,   \sin (\frac{\theta}2)   ^{8}+60\,  \sin  (\frac{\theta}2 )^{4}+540\,{\theta}^{4}-13725\,{\theta}^{2}-
7200\,\sin  ( \theta  ) \theta\\
           &&+83700\,   \sin  ( \frac{\theta}2)   ^{2} \Big) \\
L_5(\theta)&=&{\frac {\sin (\frac{\theta}2)}{1296000}}  \Big( 40500 \sin ( \theta ) {\theta}^{2}-584 \sin (\frac{\theta}2)    ^{4}\sin ( \theta ) +12000  \sin (\frac{\theta}2)    ^{4}\theta\\
&&-272   \sin (\frac{\theta}2)    ^{6}\sin (\theta ) -549000{\theta}^{3}-2745000  \sin (\frac{\theta}2)   ^{2}\theta\\
&&+44270   \sin (\frac{\theta}2)    ^{2}\sin ( \theta ) +7102095
(\theta-\sin ( \theta ))
 -64   \sin (\frac{\theta}2)   ^{8}\sin ( \theta ) +10800{\theta}^{5} \Big) 
\ee
The next formula are too large to be written here. 
We can compute also $L_{m,n}(\theta)$ for small values of $m,n$. 
For any $n\geq 2$, $\sum_{k=1}^n B_{n,k}(\theta)=1$. Since $B_{2,2}=B_2$ is know, so do $B_{2,1}$. The next ones are 
\be
L_{3,1}(\theta)&=& -\frac{1}6  \sin (\frac{\theta}2)   ^{3} \Big( 48
  \sin (\frac{\theta}2)  ^{5}\cos (\frac{\theta}2 ) -48\,  \sin (\frac{\theta}2) 
  ^{6}-12  \sin (\frac{\theta}2)   ^{3}\cos (\frac{\theta}2) +8 \sin (\frac{\theta}2)   ^{4}\\
&&+12  \sin (\frac{\theta}2)  ^{2}\theta-10\cos (\frac{\theta}2) \sin (   \frac{\theta}2 ) +40 \sin (\frac{\theta}2)  ^{2}-15\theta \Big) \\
L_{4,1}(\theta)&=&-\frac{1}{45}  \sin (\frac{\theta}2)   ^{4} \Big( 48
 \sin (\frac{\theta}2)   ^{4}\cos (\frac{\theta}2) \theta-48 \sin (\frac{\theta}2) 
  ^{5}\theta-12  \sin (\frac{\theta}2) 
  ^{5}-12\,\cos (\frac{\theta}2)   \sin (\frac{\theta}2 )   ^{2}\theta\\
&&-48 \sin (\frac{\theta}2)   ^{3}\theta+190 \sin (\frac{\theta}2)   ^{3}+9\cos (\frac{\theta}2) \theta
+96\sin (\frac{\theta}2) \theta-210\sin (\frac{\theta}2)  \Big) 
\ee
the next ones are too large again, to be written here. I am able to compute $L_{n,m}$ for $n\leq 5$ (which provides the values of $P_D^{n,m}$ for $n\leq 7$).

Using these formulae, one finds the following explicit cute values for $P_D^n$
\ben
\begin{array}{|c|c|c|c|c|c|}
\hline
n & 4 &5 &6 &7&8 \\
\hline
1-P_D^n  &\frac{35}{12\pi^2} & \frac{305}{48\pi^2} &\frac{146400\pi^2-473473}{11520\pi^4} &\frac{512400\pi^2-2900611}{23040\pi^4}&\frac {62664108221+1721664000\,{\pi }^{4}-
18670881600\,{\pi }^{2}}{48384000 \pi ^{6}}
 \\
\hline
\end{array}
\een
By Theorem \ref{theo:circle}, we can also compute $P^{4,3}_D=\frac{35}{12\pi^2}$ (or by $P_D^{4,4}+P_D^{4,3}=1$),  $P_D^{5,5}$ is given in the array, $P^{5,4}_D=\frac{65}{12\pi^2}$, $P^{5,3}_D=\frac{15}{16\pi^2}$, $P_D^{6,5}=\frac{3120\pi^2-17017}{288 \pi^4}$, $P_D^{6,4}=\frac{7200 \pi^2+57057}{3840 \pi^4}$, $P_D^{6,3}=\frac{1001}{320 \pi^4}$.

Explicit results for bi-pointed half disk ($B_1(\theta)=1$ for any $\theta\in(0,2\pi)$:
\ben
\begin{array}{|c|c|c|c|c|c|}
\hline
n &2 &3 &4 &5 &6 \\
\hline
1-\FD_n(\pi) &\frac{16}{3\pi ^2} & \frac {26}{3\pi ^2}&\frac {-83968+13725\,{\pi }^{2}}{540\pi ^{4}}&\frac {-97091+12200\,{\pi }^{2}}{240\pi ^{4}}& \frac {30749622272-4885982325\,{\pi }^{2}+
201757500\,{\pi }^{4}}{2268000{\pi }^{6}}
 \\
\hline
\end{array}
\een
The only simple formula which appears is the following:
\[\lim_{t\to 0} B_n(t)=\frac{12^n}{(n+1)(2n+1)!}.\]
which holds for all $1\leq n \leq 11$. It corresponds to the limit for $\BSEG$ with an angle going to 0.

Apart the results exposed in Theorem \ref{theo:circle}, the only explicit results in the literature concerns triangles and parallelogram (we here discuss only results known for any $n$, in 2D). 
Valtr \cite{Valtr-P} (1995) has obtained that if $S$ is a square (or a (non flat) parallelogram) then, for $n\geq 1$, 
\ben
P_S^n =\l(\frac{\binom{2n-2}{n-1}}{n!}\r)^2,
\een
and in a second paper, \cite{Valtr-T} (1996) he proved that if $T$ is a (non flat) triangle then, for $n\geq 1$,  
\ben P_T ^n= \frac{2^n (3n-3)!}{(n-1)!^3(2n)!}.
\een
Buchta \cite{BC2} goes further and gives an expression for $P_S^{n,m}$ and $P_P^{n,m}$ as a finite sum of explicit terms. 

For the bi-pointed triangle, B\'ar\'any, Rote, Steiger, Zhang \cite{BRSZ} (2000) have shown the following.  Let $T=(A,B,C)$ be a (non flat) triangle, and let $(z_1,\dots,z_n)$ be $`P_T^n$ distributed, and let $\overline{\sz[n]}=(A,B,z_1,\dots,z_n)$ be the $n+2$ tuple obtained by adding $A,B$ to $\sz[n]$. For any $n\geq 0$,
\ben\label{eq:bi1}
`P_T^n(\overline{\sz[n]}\in\CP_{n+2})= \frac{2^n}{n!(n+1)!}.
\een
These results are at the origin of numerous works concerning limit shape for convex bodies in a domain (\cite{BRSZ}, B\'ar\'any \cite{BAR}) and for the evaluation of the probability that $n$ points chosen in a convex domain $H$ are in convex position (see B\'ar\'any \cite{BAR}). 

Buchta  (2007) \cite{BC1} goes forward and prove the following result~: For any $n\geq 1$, any $1\leq m \leq n$,
\ben  \label{eq:bi2}
`P_T^n(\overline{\sz[n]}\in\CP_{n+2,m+2})=\sum_{C\in \comp(n,m)} 2^m \prod_{i=1}^m \frac{C_i}{SC_i (1+SC_{i})}
\een 
where $SC_i=C_1+\dots+C_i$ and $\comp(n,m)$ is the set of compositions of $n$ in $m$ non empty parts (Examples : $\comp(2,3)=\varnothing$, $\comp(4,2)=\{(1,3),(3,1),(2,2)\}$).

\paragraph{Additional references}
 
The literature concerning  the question of the number of points on the convex hull for i.i.d. random points taken in a convex domain is huge. I won't make a survey here but rather sends the interested reader to Reitzner \cite{Rei}, Hug \cite{Hug} and to the various paper cited in the present paper
I will focus  on what concerns the disk. 

Blaschke (1917) \cite{Bla} proves that for the 4 points problem (the so-called \it problem of Sylvester\rm ), we have for any convex $K$,
\[P_T^4\leq P_K^4\leq P_D^4.\]
B\'ar\'any (2000) \cite{BAR} have shown that 
\ben\label{eq:Bar}
\lim_{n\to +\infty} n^2(P_K^n)^{1/n}=e^2 A^3(K)/4
\een
where $A^3(K)$ is the supremum of the \it affine perimeter \rm of all convex sets $S\subset K$. 
For the disk one gets 
\ben 
 \log (P_K^n) =-2n \log n +n\log \l(2\pi^2 e^2\r)-2`e_0(3\pi^4n)^{1/5}+...
\een
where the last term, not really proved in the mathematical sense, has been obtained 
by Hilhorst \& al. \cite{HCS}. Central limit theorems exists also for the number of points on $\BCH(\sx\{n\})$ under $`P_D^n$ (and for more general domain, under the uniform or Poisson distribution), see Groeneboom \cite{Gro}, Pardon \cite{par}, B\'ar\'any and Reitzner \cite{BR}.  
 
\section{\color{red} Proof of Theorem \ref{theo:circle}}
\subsection{Proof of $(i)$} 
All along this section $n\geq 1$ is fixed. Take a closed disk $\bar{B}=\bar{B}((0,0),R_c)$, with center $(0,0)$ and radius $R_c= {1}/{\sqrt{\pi}}$, that is with area 1, and pick $n$  i.i.d. uniform points $U_1,\dots,U_n$ in $\bar{B}$. Now consider the smallest disk  $\bar{B}((0,0),R_n)$ that contains all the $U_i$'s. Clearly 
\[R_n=\inf\{ r : \#(\bar{B}(0,r)\cap \{U_1,\dots,U_n\}\}) = n\}.\]
\begin{pro}\label{pro:etape1}
Conditionally on $R_n=r$, there is a.s. exactly one index $J \in \cro{1,n}$ such that $U_J$ belongs to the circle ${B}((0,0),r)$. Conditionally on $\{J=j, R_n=r\}$, $U_j$ and $(U_1,\dots,U_{j-1},U_{j+1},\dots,U_n)$ are independent, $U_j$ has the uniform law on the circle $\partial B((0,0),r)$, and $U_1,\dots,U_{j-1},U_{j+1},\dots,U_n$ are uniform in $\bar{B}((0,0),r)$. 
\end{pro}
\Proof A.s. the points $U_1,\dots,U_n$ are not on the same circle with center $(0,0)$, and by symmetry conditionally on $R_n=r$ and $J=j$,  $U_j$ is uniform on $B((0,0),r)$. Now, conditionally on $R_n=r$ and $J=j$, each variable $U_\ell$ (for $\ell \neq j$) are just conditioned to satisfy $\|U_\ell\|_2\leq r$, and this conditioning  conserves the uniform distribution. ~$\Box$ \medskip

\noindent {\sl Proof of Theorem \ref{theo:circle}$(i)$}. Theorem \ref{theo:circle}$(i)$ is -- or should be -- intuitively obvious, taking into account Proposition \ref{pro:etape1}. This Proposition says that the two following models $(a)$ and $(b)$:\\
-- $(a)$  $n$ points i.i.d. uniform in a disk,\\
-- $(b)$  one point uniform on the circle and, independently, $n-1$ i.i.d. uniform inside the disk\\
are equivalent with respect to the probability to be in convex position. \\
Now if we come back to the $\BSEG$ considerations, when $\theta\to 2\pi$,  the points $w_1(R_c,\theta)$ and $w_2(R_c,\theta)$ become closer and closer, and the line passing by these points lets all the other points in one of the half plane it defines. It is intuitively clear that replacing $w_1(R_c,\theta)$ and $w_2(R_c,\theta)$ by a single point \it close \rm to them (for example,  at position $(-R_c,0)$) will not dramatically change the model nor the probability to be in convex position. This is the essence of Theorem \ref{theo:circle}$(i)$. \par
For sake of completeness, let us give a formal proof. Take $R>0$ and consider the two sets $S(`e)=\SEG(2\pi-`e,R)$ and $S=\SEG(2\pi,R)=\bar{B}((0,0),R)$. 
These two sets are close for the Hausdorff topology when $`e$ is small. We always have $S(`e)\subset S$, and $|S\setminus S(`e)|$ goes to 0. This property 
implies that if we fix  $`e'>0$, for $`e$ small enough, for $z_1,\dots,z_n$ chosen uniformly and independently under $`P_{S}$, 
 \ben\label{eq:lambda}
 `P(\{z_1,\dots,z_n\}\subset S(`e))\geq 1-`e'. 
 \een
Conditionally on the event  $\Lambda_{`e}:=\big\{\{z_1,\dots,z_n\}\subset S(`e)\big\}$, the $z_i$'s are i.i.d. uniform in $S(`e)$.
Let $w_1(`e)=w_1(R_c,2\pi-`e)$, $w_2(`e)=w_2(R_c,2\pi-`e)$, $w=-R$.  

We want to show that $`P((z_1,\dots,z_n,w_1^{`e},w_2^{`e})\in \CP_{n+2}|\Lambda_{`e})\to`P((z_1,\dots,z_n,w)\in \CP_{n+1})$.
Consider the following sets (subsets of $S^n$):
\be
E_1(`e)&:=&\{(t_1,\dots,t_n) \in S(`e)~: (t_1,\dots,t_n,w_1(`e),w_2(`e)) \in \CP_{n+2}\}\\
E_2&:=&\{(t_1,\dots,t_n) \in S~:  (t_1,\dots,t_n,w)\in \CP_{n+1}\}.
\ee
It suffices to prove that $|E_1(`e)|\sous{\to}{`e \to 0} |E_2|$.
First $E_1(`e)\subset E_2$ since if $(t_1,\dots,t_n,w_1(`e),w_2(`e))$ belongs to $\CP_{n+2}$ and since the segments $[w_1(`e),w]$ and $[w_2(`e),w]$ are chords, then $ (z_1,\dots,z_n,w_1(`e),w_2(`e),w)$ is in $\CP_{n+3}$ from what we deduce that $E_2$ is in $\CP_{n+1}$. \par
To end the proof take $(t_1,\dots,t_n)\in E_2$. We show that when $`e$ is small enough, it is in $E_1(`e)$. 
More precisely, we will see that it is not the case only if the $t_i$ belongs to a null set (for Lebesgue measure).
We assume that $n\geq 2$ since for $n=1$ the result is clear. 

First, for $`e>0$ small enough, if the $t_i$'s are different and different to $-R$, all the $t_i$ belongs to $S(`e)$. 
Since $(t_1,\dots,t_n,w) \in \CP_{n+1}$ draw the convex polygon $p$ passing by these points, and relabel the $t_i'$s as $t_1^\star,\dots,t_n^\star$ clockwise around $p$ so that  the neighbours of $w$ are  $t^\star_1$ and $t^\star_n$. Again, up to null set, the angles $(w,t^\star_1,t^\star_2)$ and $(t^\star_{n-1},t^\star_n,w)$ are not 0, and it appears clearly that for $`e$ small enough, $(t_1,\dots,t_n,w_1(`e),w_2(`e))\in \CP_{n+2}$.
We then have 
$E_2=\cup_{`e} E_1(`e)$ and the $E_1(x)\cup E_1(x')$ if $x'<x$, so $|E_1(`e)|\to |E_2|$ when $`e$ goes to 0. ~$\Box$.

\subsection{Proof of  $(ii)$}

For any $\theta\in[0,2\pi], R>0$,
\ben\label{eq:seg}
|\SEG(\theta,R)|:=\frac{R^2}2\l(\theta-\sin(\theta)\r)
\een
and then for 
\ben\label{eq:Rtheta}
R_\theta=\sqrt{\frac{2}{\theta-\sin(\theta)}},
\een
the area $|\SEG(\theta,R_\theta)|=1$. Denote more simply by $\SEG_\theta$ the segment $\SEG(\theta,R_\theta)$ with unit area. The size $L_\theta$ of the special border $[w_1(\theta,R_\theta),w_2(\theta,R_\theta)]$ for this segment is
\ben\label{eq:Ltheta}
L_\theta=2 R_\theta \sin(\theta/2).
\een

In this section we fix $\theta\in(0,2\pi)$ and search to express $B_n(\theta)$ with some combinations of $B_j(\nu)$, for $\nu< \theta$ and $j< n$. To get the decomposition we will ``push the arc of circle'' $\AC(\theta,R)$ inside $\SEG(\theta,R_\theta)$ till it touches one of the $z_i$'s doing something similar to the Buchta's method (for the computation of $P_S^n$ and $P_T^n$). Here it is a bit more complex: we need the arc of circle to stay an arc of circle during the operation in order to get a nice decomposition, and also we somehow need to keep the bi-pointed elements. The arc angle and radius will change during the operation. This will lead to a quadratic formula for $B_n$.
\begin{figure}[ht]
\centerline{\includegraphics{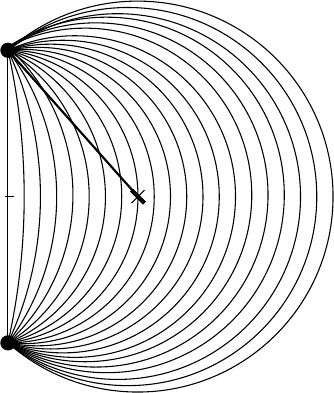} ~~\includegraphics {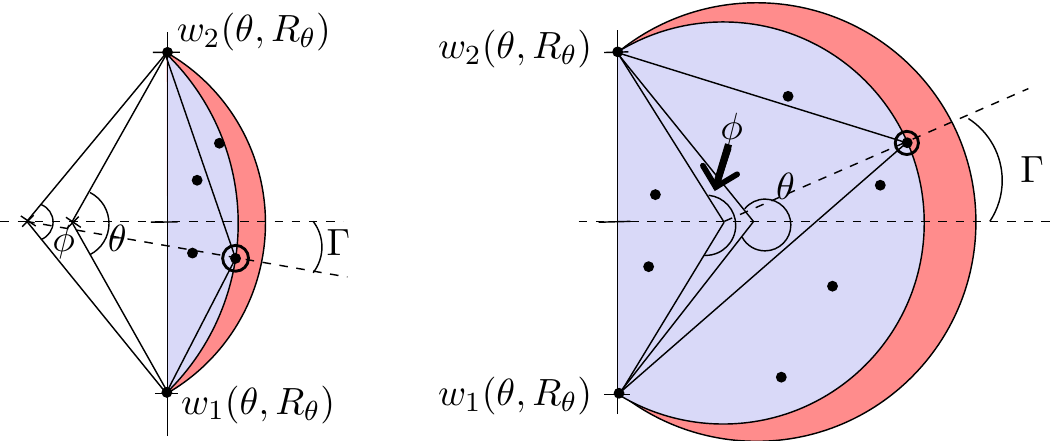}}
\captionn{\label{fig:filled-3} Representation of the family ${\cal F}_\theta$. The angle $\phi<\theta$ and $\SEG[\phi]\leq \SEG[\theta]$. The angles are taken at the center of the circle that defines the segments. }
\end{figure}
Almost of quantities appearing in this section should be indexed by $\theta$.  In order to avoid heavy notation we won't do this. Draw $\SEG_\theta$ in the plane. We consider the family of segments 
\[{\cal F}_\theta:=(\SEG[\phi],0\leq \phi\leq \theta)\]
 having as special border the special border of $\SEG_\theta$, that is $[w_1(\theta,R_\theta),w_2(\theta,R_\theta)]$, and lying at its right, such that the angle of $\SEG[\phi]$ is $\phi$ (see Fig. \ref{fig:filled-3}). \par
When $\phi$ goes from $\theta$ to $0$, the center $O[\phi]$ of (the circle which defines) $\SEG[\phi]$ moves on the $x$-axis from $O[\theta]=0$ to $(-\infty,0)$. Comparing the distance from $O[\phi]$ to the special border, we can compute the coordinate of $O[\phi]$:
\ben\label{eq:center}
O[{\phi}]=\frac{L_{\theta}}2 \l(\cot(\frac{\theta}{2})-\cot(\frac{\phi}{2})\r)
\een
and the radius of $\SEG[\phi]$,
\ben\label{eq:Rphi} R[\phi]=R_\theta\frac{\sin(\theta/2)}{\sin(\phi/2)}.\een
Since the special border of all the $\SEG[\phi]$ is the same one sees that if $\phi <\phi'$ then $\SEG[\phi]\subset \SEG[\phi']$. When $\phi$ goes to 0, $\SEG[\phi]$ goes to $[w_1(\theta,R_\theta),w_2(\theta,R_\theta)]$ (for the Hausdorff topology). One also sees that $\SEG[\theta]=\SEG_\theta$, and for $\phi<\theta$, by  \eref{eq:seg} and \eref{eq:Rtheta} , 
\ben\label{eq:phit}
|\SEG[\phi]|&=&\l(\frac{\sin(\theta/2)}{\sin(\phi/2)}\r)^2\frac{\phi-\sin(\phi)}{ \theta-\sin(\theta)} 
\een
and then the other segments of the family ${\cal F}_\theta$ have area smaller than 1 (see Fig. \ref{fig:filled-3}). 

Again $\theta$ is fixed. Let $z_1,\dots,z_n$ be $n\geq 1$ i.i.d. uniform random points in $\SEG_\theta$.
Denote by 
\[\Phi=\min\l\{\phi: \#(\{z_1,\dots,z_n\} \cap \SEG[\phi]\r\})=n\},\]
and let $J$ the (a.s. unique) index of the variable $z_j$ on $\partial \SEG[\phi]$.
Finally let $\Gamma$ be the (signed) angle $\big((+\infty,0),O[\Phi],z_J\big)$  formed by the $x$-axis and the line $(0[\Theta],z_J)$ (see Fig. \ref{fig:filled-3}). We have
\begin{pro}\label{pro:dens} The distribution of $(\Phi,\Gamma)$ admits the following density $f_{(\Phi,\Gamma)}$ with respect to the Lebesgue measure 
 \[f_{(\Phi,\Gamma)}(\phi,\gamma)= n\, 
\frac{\sin(\theta/2)^{2n}}
     {(\theta-\sin(\theta))^{n} } 
\frac{(\phi-\sin(\phi))^{n-1}}
     {\sin(\phi/2)^{2n+1}}
\l(\cos(\gamma)-\cos(\phi/2)\r)  1_{0\leq \phi\leq \theta} \,1_{|\gamma|\leq \phi/2}     .\]
\end{pro}

\Proof First, the density of $z_J=(x,y)$ with respect to the Lebesgue measure on $|\SEG_\theta|$ is $n dxdy |\SEG_{x,y}|^{n-1}$ where $|\SEG_{x,y}|^{n-1}$ is the area of the unique element of the family ${\cal F}_\theta$ whose border contains $(x,y)$. We then just have to make a change of variables in this formula ! 

We search the unique pair $(\phi,\gamma)$ such that 
\[x+iy=R[\phi]e^{i\gamma}+O[\phi].\]
Since by \eref{eq:phit} and \eref{eq:center} everything is explicit, we can compute the Jacobian 
\[\l| {\sf det}\l(
\begin{array}{cc}
\frac{\partial x}{\partial \phi} & \frac{\partial x}{\partial \gamma}\\   
\frac{\partial y}{\partial \phi} & \frac{\partial y}{\partial \gamma}   
\end{array}\r) \r| =\frac{\sin(\theta/2)^2}{\sin(\phi/2)^3}\frac{(\cos(\gamma)-\cos(\phi/2))}{(\theta-\sin(\theta))}.
\]
From what we deduce the wanted formula, using \eref{eq:phit}.~$\Box$.
 
Now, it remains to end the decomposition of our problem.  Conditionally on  $(\Phi,\Gamma,J)=(\phi,\gamma,j)$, 
the points $z_1,\dots,z_{j-1},z_{j+1},\dots,z_{n}$ are i.i.d. uniform in $\SEG[\phi]$.

\begin{figure}[ht]
\centerline{\includegraphics{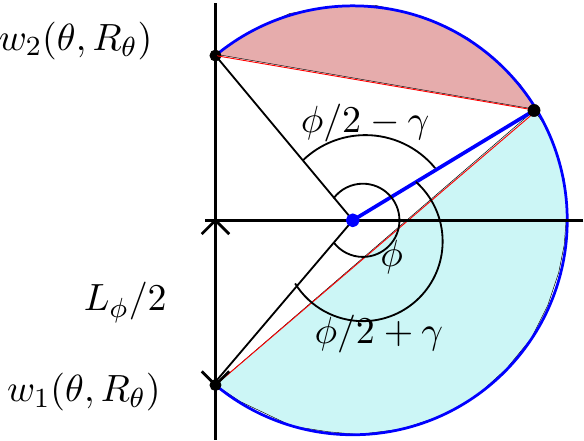}}
\captionn{\label{fig:bi-segments}Decomposition of the computation of $B_n(\theta)$, and definition of the two sub-segments appearing in the decomposition.
  }
\end{figure}

The triangle $T:=(w_1(\theta,R_\theta),w_2(\theta,R_\theta),z_j)$ is inscribed in $\SEG[\phi]$ and $\SEG[\phi]\setminus T$ produces 
two segments $S_1$ and $S_2$. Since we may rescale $\SEG[\phi]$ to be $\SEG_\phi$ (to get area 1), the question now is that of the area of the two rescaled segments. 
After rescaling, $S_1$ and $S_1$ appear to be  $\SEG[\phi /2+\gamma,R_\phi ]$ and $\SEG[\phi /2+\gamma,R_\phi ]$ since these lards have the right angles. Using \eref{eq:seg}
\ben\label{eq:phi-alpha}
|\SEG[ \alpha,R_\phi]|&=&    \frac{\alpha-\sin(\alpha)}{\phi-\sin(\phi)}.
\een
We keep temporally  notation $S_1$ and $S_2$ instead of $\SEG[\phi /2+\gamma,R_\phi ]$ and $\SEG[\phi /2+\gamma,R_\phi ]$ for short.
The following Proposition is a simple consequence of the fact that the uniform distribution is conserved by conditioning. It is the ``combinatorial decomposition'' of the computation of $B_n(\theta)$, illustrated on Fig. \ref{fig:bi-segments}. 
\begin{pro}\label{pro:er}
$(i)$ Conditionally on $(\Phi,\Gamma,J)=(\phi,\gamma,j)$, 
the respective number $(N_1,N_2,N_3)$ of points of $z\{n\}\setminus\{z_j\}$ in $S_1$, $S_2$ and $\SEG_\phi-(S_1\cup S_2)$   
 is 
\[{\sf Multinomial}(n-1,|S_1|,|S_2|,1-|S_1|-|S_2|).\]
$(ii)$ Conditionally on $(\Phi,\Gamma,J)=(\phi,\gamma,j)$ and $(N_1,N_2,N_3)=(k_1,k_2,k_3)$ the points $z_1,\dots,z_n$ are in convex position with probability 
$1_{k_3=0,k_1+k_2=n-1}B_{k_1}(\phi /2+\gamma)B_{k_2}(\phi/2-\gamma)$.
\end{pro}

Putting everything together we have obtained
\be
B_n(\theta)&=&\int_{0}^{\theta}\int_{-\frac{\phi}2}^{\frac{\phi}2}f_{(\Phi,\Gamma)}(\phi,\gamma)
\sum_{k=0}^{n-1}\binom{n-1}k
|S_1|^k|S_2|^{n-1-k}
B_k(\phi/2+\gamma)
B_{n-1-k}(\phi/2-\gamma) d\gamma d\phi
\ee
Set $\eta=\phi/2+\gamma$, $d\eta=d\gamma$, $\eta$ goes from $0$ to $\phi$ (and $\phi/2-\gamma=\phi-\eta$),
giving
 \ben 
B_n(\theta)&=&\int_{0}^{\theta}\int_{0}^{\phi}f_{(\Phi,\Gamma)}(\phi,\eta-\phi/2)
\sum_{k=0}^{n-1}\binom{n-1}k\\
&\times & |\SEG[\eta,R_\phi]|^k|\SEG[\phi -\eta,R_\phi]|^{n-1-k}
B_k( \eta)
B_{n-1-k}(\phi-\eta)d\eta d\phi
\een
from which we get
\ben\label{eq:F_n}
B_n(\theta)&=&\int_{0}^{\theta}\int_{0}^{\phi} 
 n 
\frac{\sin(\theta/2)^{2n}}
     {(\theta-\sin(\theta))^{n} } 
\frac{\cos(\eta-\phi/2)-\cos(\phi/2)}
     {\sin(\phi/2)^{2n+1}}\sum_{k=0}^{n-1}\binom{n-1}k\\
&\times& 
  (\eta-\sin(\eta))^kB_k( \eta)((\phi-\eta)-\sin(\phi-\eta))^{n-1-k}B_{n-1-k}(\phi-\eta) d\eta d\phi.
\een
Now, $\cos(\eta-\phi/2)-\cos(\phi/2)=2\sin(\eta/2)\sin((\phi-\eta)/2)$.
Finally setting $L_n(\theta)$ as done in \eref{eq:F-L}, we obtain
Formula \ref{theo:circle}$(ii)$. 
 
\subsection{Proof  $(i')$}
\label{eq:iprime}

Recall Proposition \ref{pro:etape1}. To compute $P_n^D$ we can work under the model where $n-1$ points $z_1,\dots,z_{n-1}$ 
are picked independently and uniformly  
inside the disk $B((0,0),R_c)$  (with $R_c=\pi^{-1/2}$) and one point on the border. We place this last point at position $-R_c$ which is allowed since rotation keeps convex bodies and the uniform distribution.  

\begin{figure}[ht]
\centerline{\includegraphics{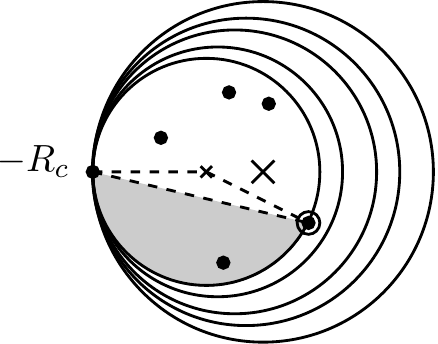}}
\captionn{\label{fig:czd}Decomposition of the computation of $P_n^D$. The big cross is the center of the initial circle, the small one, the center of the smallest circle containing all the points.  }
\end{figure}

Now take a family of circles ${\cal G}=\{B[r],0\leq r \leq R_c\}$ such that $B[r]$ as radius $r$, its center at position $-R_c+r$, implying that $-R_c$ belongs to all these circles (see Fig. \ref{fig:czd}).

If $r'<r$, $B[r']\subset B[r]$. Let $r^\star$ be the largest 
circle such that exists $1\leq k\leq n-1$, $z_k\in \partial B[r^\star]$. Denote then by $\phi$ the angle such that $z_k=(-R_c+r)+re^{i(-\pi+\phi)}$.
If we denote by $(X,Y)$ the (Euclidean) position of $z_k$, the density of the distribution of $(x,y)$ is 
\[(n-1)1_{(x,y)\in B((0,0),R_c)} |B[r]|^{n-2}    dxdy \] 
where $B[r]$ is the unique circle in the family ${\cal G}$ which passes by $(x,y)$. 
We can then compute the Jacobian and find the distribution of 
$(r,\phi)$ to be with density $1_{0\leq r \leq R_c, 0\leq \phi\leq 2\pi} r(1-\cos(\phi)) (\pi r^2)^{n-2} dr d\phi$. 
Once $z_k$ is given, we can once again normalise the problem, and come back on a circle of area $R_c$. 
We then get, using $1+\cos(\phi)=2\sin^2(\phi/2)$
\be
P_D^{n} &=&(n-1) \int_{0}^{R_c} \int_{0}^{2\pi}\sum_{k=0}^{n-2}\binom{n-2}k 2\sin^2(\phi/2)  r(\pi r^2)^{n-2}\\
&\times& B_k(\phi)B_{n-2-k}(2\pi-\phi) |\SEG(\phi,R_c)|^k |\SEG(2\pi-\phi,R_c)|^{n-2-k} d\phi dr
\ee
The integration with respect to $dr$ gives
\be
P_D^{n} &=&\frac{1}{\pi} \int_{0}^{2\pi}\sum_{k=0}^{n-2}\binom{n-2}k  \sin^2(\phi/2)
 B_k(\phi)B_{n-2-k}(2\pi-\phi) \\
&\times &\l(\frac{\phi-\sin(\phi)}{2\pi}\r)^k \l(\frac{2\pi-\phi+\sin(\phi)}{2\pi}\r)^{n-2-k}   d\phi 
\ee
since once $\phi$ is known, the convexity follows that on the pair of bi-pointed segments 
with angles $\phi$ and $2\pi-\phi$, and the number of elements in these segments is ${\sf binomial}\l(n-2,|\SEG(\phi,R_c)|\r)$. 
 
\subsection{Proof of $(iii)$}
The proof is the same as that of $(ii)$ except that in Proposition \ref{pro:er} we need to follow the number of points falling in the triangle. We then get  
\be 
B_{n,m}(\theta)&=&\int_{0}^{\theta}\int_{0}^{\phi}f_{(\Phi,\Gamma)}(\phi,\eta-\phi/2)
\sum_{n_1+n_2+n_3=n-1\atop{m_1+m_2=m-1}}
\binom{n-1}{n_1,n_2,n_3}\\
&\times & |\SEG[\eta,R_\phi]|^{n_1}|\SEG[\phi -\eta,R_\phi]|^{n_2}(1-|\SEG[\eta,R_\phi]|-|\SEG[\phi -\eta,R_\phi]|)^{n_3}\\
&\times&
B_{n_1,m_1}( \eta)
B_{n_2,m_2}(\phi-\eta)d\eta d\phi
\ee
Using the notation introduced in \eref{eq:F-L} we get $(iii)$ .

\subsection{Proof of $(iii')$}
Copy the arguments in Section \ref{eq:iprime}.
In the same way, for $n\geq 2$, $1\leq m\leq n$ 
\be
P_D^{n,m}&=&\frac{1}{\pi} \int_{0}^{2\pi}\sum_{k=0}^{n-2}\sum_{1\leq m_1 \leq n-2}\binom{n-2}k \sin^2(\phi/2)  \\
&\times& B_{k,m_1}(\phi)B_{n-2-k,m-m_1-2}(2\pi-\phi) \l(\frac{\phi-\sin(\phi)}{2\pi}\r)^k \l(\frac{2\pi-\phi+\sin(\phi)}{2\pi}\r)^{n-2-k}   d\phi
\ee
with the condition that $B_{k,k+l}=0$.    $(iii')$ follows.

\subsection{Proof of $(iv)$}
The same proof of $(i)$ does the job.

\section{Effective computation of $L_n$}
\label{sec:CC}
We explain in this part how to make effective computations. 
Since $B_1(\theta)=B_0(\theta)=1$, $L_0(\theta)$ and $L_1(\theta)$ are known by \eref{eq:F-L}.  
Bruno Salvy \cite{BS}  in a personal communication gave me a method to compute $B_n$ (my personal method fails at $n=7$).

Denote by $J_n(t)=\int_{0}^{t}  \sum_{k=0}^{n-1}B_{k}(u )B_{n-1-k}(t-u)    du$,
and by $TJ_n$, $TB_n$ the Laplace transform of $J_n$ and $B_n$. 
We have 
\[TJ_n(s)=\sum_{k=0}^{n-1} TB_k(s) TB_{n-k-1}(s).\]
It turns out that knowing the first values of $B_k$, the computation of $J_n$ by the previous formula and by inversion of the Laplace transform ends, using maple (when it does not by simple integration).
Then it appears that $J_n(v)$ is a polynomial in $\sin(v), \cos(v)$ and $v$. The subsequent integration $\int_{0}^t J_n(v)/\sin(v)^{2n+1}dv$ is possible by helping the computer. Using $\cos(v)^2+\sin(v)^2=1$, it is possible to rewrite $J_n$ as  polynomial of degree at most 1 in $\cos(v)$. Then write $J_n(v)/\sin(v)^{2n+1}$ under the form $Q_0+Q_1/\sin(v)^1+\sum_{l=1}^{2n+1} Q_l/\sin(v)^l$, where $Q_n$ is a polynomial in $\cos(v)$ and  $v$. Then, proceed to successive integrations by parts, starting from the largest degree at the denominator till $l=2$. The remaining integral is computed by a simple integration.

This method allows one to compute 11 terms with maple. 
The computation of $B_{n,m}$ is possible using the same algorithm, except that some complications arise from the inverse Laplace which makes appear some polylogarithm functions (of the type ${\sf polylog}(n,e^{it})+{\sf polylog}(n,-e^{it}))$ for some $n\geq 2$). There limit at $t\to 0$ have to be treated separately, letting $\zeta(n)$ plays a role in the final result.  I am able to compute $B_{n,m}$ for $n\leq 5$, and then $`P_{n,m}^D$ for $n\leq 7$.

\small
\bibliographystyle{abbrv}
\bibliography{bibfile}

\begin{thebibliography}{10}

\bibitem{BAR}
I.~B\'ar\'any.
\newblock Sylvester's question: The probability that n points are in convex
  position.
\newblock {\em Ann. Probab.}, 27(4):2020--2034, 1999.

\bibitem{BR}
I.~B\'{a}r\'{a}ny and M.~Reitzner.
\newblock Poisson polytopes.
\newblock {\em The Annals of Probability}, 38(4):1507--1531, 07 2010.

\bibitem{BRSZ}
I.~B\'ar\'any, G.~Rote, W.~Steiger, and C.~Zhang.
\newblock A central limit theorem for random convex chains.
\newblock {\em Discrete comput. Geom}, 30:35--50, 2000.

\bibitem{Bla}
W.~Blaschke.
\newblock \"{U}ber affine geometrie xi: L\"{o}sung des 'vierpunktproblems' von
  sylvester aus der theorie der geometrischen wahrscheinlichkeiten.
\newblock {\em Ber. Verh. Sachs. Akad. Wiss. Leipzig Math.-Phys}, 69:436--453,
  1917.

\bibitem{BC1}
C.~Buchta.
\newblock The exact distribution of the number of vertices of a random convex
  chain.
\newblock {\em Mathematika}, 53(2):247--254 (2007), 2006.

\bibitem{BC2}
C.~Buchta.
\newblock On the number of vertices of the convex hull of random points in a
  square and a triangle.
\newblock {\em Anz. \"Osterreich. Akad. Wiss. Math.-Natur. Kl.}, 143:3--10,
  2009/10.

\bibitem{Gro}
P.~Groeneboom.
\newblock Limit theorems for convex hulls.
\newblock {\em Probability Theory and Related Fields}, 79(3):327--368, 1988.

\bibitem{HCS}
H.~Hilhorst, P.~Calka, and G.~Schehr.
\newblock Sylvester's question and the random acceleration process.
\newblock {\em Journal of Statistical Mechanics: Theory and Experiment},
  2008(10):P10010, 2008.

\bibitem{Hug}
D.~Hug.
\newblock Random polytopes.
\newblock In {\em Stochastic geometry, spatial statistics and random fields},
  volume 2068 of {\em Lecture Notes in Math.}, pages 205--238. Springer,
  Heidelberg, 2013.

\bibitem{par}
J.~Pardon.
\newblock Central limit theorems for uniform model random polygons.
\newblock {\em Journal of Theoretical Probability}, 25(3):823--833, 2012.

\bibitem{Rei}
M.~Reitzner.
\newblock {\em New perspectives in stochastic geometry, Molchanov, I., and
  Kendall, W. edts}, chapter Random polytopes (survey)., pages 45--76.
\newblock Oxford University Press, Oxford., 2010.

\bibitem{BS}
B.~Salvy.
\newblock How to compute $l_n$.
\newblock Personal communication, 2013.

\bibitem{Valtr-P}
P.~Valtr.
\newblock Probability that n random points are in convex position.
\newblock {\em Discrete {\&} Computational Geometry}, 13:637--643, 1995.

\bibitem{Valtr-T}
P.~Valtr.
\newblock The probability that {\it n} random points in a triangle are in
  convex position.
\newblock {\em Combinatorica}, 16(4):567--573, 1996.

\end{thebibliography}
 \end{document}